\documentclass[reqno, a4paper, 10pt]{amsart}
\usepackage{amssymb,url, color, mathrsfs}
\usepackage[colorlinks=true, bookmarks=true, pdfstartview=FitH, pagebackref=true, linktocpage=true]{hyperref}
\usepackage[short,nodayofweek]{datetime}
\usepackage[pagewise,running,mathlines,displaymath, switch]{lineno}
\usepackage{times}
\usepackage{indentfirst, tabularx}
\usepackage[table]{xcolor}
\usepackage{float}

\theoremstyle{plain}
\numberwithin{equation}{section}
\newtheorem{theorem}{Theorem}

\newtheorem{proposition}{Proposition}
\theoremstyle{remark}

\DeclareMathOperator{\Rset}{\mathbf{R}}

\definecolor{brown}{rgb}{0.5,0,0}
\definecolor{backgroundcolor}{rgb}{0.98, 0.92, 0.73}

\newif\ifprint
\ifprint
\else
\fi

\author[T.V. Duoc]{Trinh Viet Duoc}
\address[T.V. Duoc]{Department of Mathematics\\
College of Science, Vi\^{e}t Nam National University\\
H\`{a} N\^{o}i, Vi\^{e}t Nam.}
\email{\href{mailto: T.V. Duoc <tvduoc@gmail.com>}{tvduoc@gmail.com}}

\author[Q.A. Ng\^o]{Qu\^{o}\hspace{-0.5ex}\llap{\raise 1ex\hbox{\'{}}}\hspace{0.5ex}c Anh Ng\^o}
\address[Q.A. Ng\^o]{Department of Mathematics\\
College of Science, Vi\^{e}t Nam National University\\
H\`{a} N\^{o}i, Vi\^{e}t Nam.}
\email{\href{mailto: Q.A. Ng\^o <nqanh@vnu.edu.vn>}{nqanh@vnu.edu.vn}}
\email{\href{mailto: Q.A. Ng\^o <bookworm\_vn@yahoo.com>}{bookworm\_vn@yahoo.com}}

\begin{document} 

\allowdisplaybreaks

\setpagewiselinenumbers
\setlength\linenumbersep{100pt}

\title[A note on radial solutions of $\Delta^2 u + u^{-q} = 0$ in $\mathbb R^3$ with quadratic growth]
{A note on radial solutions of $\Delta^2 u + u^{-q} = 0$ in $\mathbb R^3$ with exactly quadratic growth at infinity}

\begin{abstract}
Of interest in this note is the following geometric interesting equation $\Delta^2 u + u^{-q} = 0$ in $\mathbb R^3$. It was found by Choi--Xu (J. Differential Equations \textbf{246}, 216--234) and McKenna--Reichel (Electron. J. Differential Equations \textbf{37} (2003)) that the condition $q>1$ is necessary and any radially symmetric solution grows at least linearly and at most quadratically at infinity for any $q>1$. In addition, when $q>3$ any radially symmetric solution is either exactly linear growth or exactly quadratic growth at infinity. Recently, Guerra (J. Differential Equations \textbf{253}, 3147--3157) has shown that the equation always admits a unique radially symmetric solution of exactly given linear growth at infinity for any $q>3$ which is also necessary. In this note, by using the phase-space analysis, we show the existence of infinitely many radially symmetric solutions of exactly given quadratic growth at infinity for any $q>1$.
\end{abstract}

\date{\bf \today \; at \, \currenttime}



\maketitle

\section{Introduction}

In this note, we are interested in entire solutions of the following geometric interesting equation 
\begin{equation}\label{eqMain}
\Delta^2 u + u^{-q} = 0
\end{equation}
in $\Rset^3$ with $q>0$. Recently, equations of the type \eqref{eqMain} have been captured much attention since they are natually arised when studying the prescribed Q-curvature problem either in $\mathbb R^3$ (with a flat background metric) or in $\mathbb S^3$. To be precise, positive smooth solutions of Eq. \eqref{eqMain} for the case $q=7$ correspond to conformal metrics conformally equivalent to the flat metric which have constant Q-curvature in $\mathbb R^3$. Moreover, upon using the stereographic projection, any conformal metric on $\mathbb S^3$ is simply a suitable pullback of the standard one $g_{\mathbb S^3}$ under the conformal transformation of $\mathbb S^3$ into itself; see \cite{ChoiXu}. For interested readers, we refer to \cite{ChoiXu} and the references therein.

As far as we know, Eq. \eqref{eqMain} was first studied by Choi and Xu in an preprint in 1999, which is eventually published in \cite{ChoiXu}, by Xu in \cite{xu2005}, and then by McKenna and Reichel for $\mathbb R^n$ for arbitrary $n \geqslant 3$ in \cite{KR}. To seek for complete conformal metrics on $\mathbb S^3$, it is often to look for $C^4$ positive solutions $u$ of Eq. \eqref{eqMain} with exactly linear growth at infinity in the sense that $\lim_{|x| \to +\infty}  u(x)/|x| =\alpha$ for some non-negative constant $\alpha$ in the case $q=7$. In this scenario, it is worth noticing that, $C^4$ positive solutions of Eq. \eqref{eqMain} with exactly linear growth at infinity is completely classified. Indeed, it was found by Choi and Xu that, up to a constant multiple, translation and dilation, there holds $u(x)=\sqrt{1+|x|^2}$. Then, it is natural to study $C^4$ positive solution $u$ of Eq. \eqref{eqMain} when $q \ne 7$ and when $u$ is no longer of linear growth at infinity which corresponds to incomplete conformal metrics on $\mathbb S^3$.

As a first step toward answering this question, we first look for radially symmetric solutions of Eq. \eqref{eqMain}. However, in order to understand the motivation of writing this note, we first collect all results found in \cite{ChoiXu} and in \cite{KR}. The following result is now well-known.

\begin{theorem}[see  \cite{ChoiXu, KR}]\label{thmChoiXu}
We have the following claims:
\begin{itemize}
  \item[(a)] If Eq. \eqref{eqMain} admits a smooth positive solution on $\mathbb R^3$, then there must hold $q>1$.
  \item[(b)] If Eq. \eqref{eqMain} admits a smooth positive solution on $\mathbb R^3$ with exactly linear growth, then $q>3$.
  \item[(c)] Any radially symmetric solution of Eq. \eqref{eqMain} grows \textbf{at least linearly} at infinity in the sense that $\liminf_{|x| \to +\infty}  u(x)/|x|>0$ and \textbf{at most quadratically} at infinity in the sense that $\limsup_{|x| \to +\infty}  u(x)/|x|^{2+\varepsilon} = 0$ for arbitrary $\varepsilon>0$.
  \item[(d)] If $1<q<3$, then Eq. \eqref{eqMain} admits infinitely many radially symmetric, \textbf{singular} solution with growth rate \textbf{strictly between linear and quadratic}, these are of the form $b r^\beta$ with $\beta \in (1,2)$.
  \item[(e)] If $q>1$, then there exist radially symmetric and smooth solutions of Eq. \eqref{eqMain} which grow \textbf{super-linearly} at infinity.
  \item[(f)] If $q>3$, then any radially symmetric and smooth solution of Eq. \eqref{eqMain} is either \textbf{exactly linear growth or exactly quadratic growth} at infinity.
  \item[(g)] If $q \geqslant 7$, then there exist a unique radially symmetric and smooth solutions of Eq. \eqref{eqMain} with \textbf{linear growth} at infinity.
\end{itemize}
\end{theorem}


Recently, by using the phase-space analysis, Guerra \cite{G} studied the structure of radially symmetric solutions of Eq. \eqref{eqMain} without assuming $q=7$. As far as we know, he first showed, among others, that Eq. \eqref{eqMain} also admits solutions with exactly linear growth at infinity for any $q>3$; see also \cite{Lai14} for another proof based on the variation of parameters formula for ODEs. The following is his finding.

\begin{theorem}[see \cite{G}]\label{thmGuerra}
We have the following cases:
\begin{itemize}
  \item[(a)] For $q>3$, there exists a unique radially symmetric solution of Eq. \eqref{eqMain} such that $\lim_{|x| \to +\infty}  u(x)/|x|$ exists.
  \item[(b)] For $q=3$, there exists a unique radially symmetric solution of Eq. \eqref{eqMain} such that $\lim_{|x| \to +\infty}  u(x)/ \big( |x| (\log|x|)^{1/4} \big)= 2^{1/4}$.
  \item[(c)] For $1<q<3$, there exists a unique radially symmetric solution of Eq. \eqref{eqMain} such that $\lim_{|x| \to +\infty}  u(x)/ |x|^\tau= K_q^{-1/(q+1)}$ where $K_q=\tau (2 - \tau) (\tau+1)(\tau-1)$ and $\tau = 4/(q+1)$.
\end{itemize}
\end{theorem}

In view of Theorem \ref{thmChoiXu}(f) and Theorem \ref{thmGuerra} above, the present note has twofold. First we improve Choi--Xu's result by showing that there exist radially symmetric solutions of Eq. \eqref{eqMain} with exactly quadratic growth at infinity. Second, we prove that the quadratic growth can be arbitrary. To be precise, we shall prove the following result.

\begin{theorem}\label{thmMAIN}
Given any $\kappa >0$ and any $q>1$, there exist infinitely many radially symmetric solutions $u$ of exactly quadratic growth at infinity in the sense that
$$\lim_{|x| \to +\infty}  \frac{u(x)}{|x|^2} =\kappa$$ 
holds. Furthermore, the solution $u$ and the given limit $\kappa$ are related through the following identity
\[
6\kappa = (\Delta u)(0) -\int_0^{+\infty}t u^{-q}(t) dt .
\]
\end{theorem}

Clearly, an easy consequence of Theorem \ref{thmMAIN} is that the geometric interesting equation $\Delta^2 u + u^{-7} = 0$ in $\mathbb R^3$ and its corresponding integral equation $u(x)=\int_{\mathbb R^3} |x-y| u(y)^{-7} dy$ are not equivalent since the latter equation only admits radial solutions with linear growth at infinity; see \cite[Theorem 1.1]{xu2005}. Further investigation for the relation between these two equations will be carried out in future. In addition, Theorem \ref{thmMAIN} shows that at infinity, the highest order term of $u$ is $|x|^2$. In the next result, we study lower order terms of $u$ at infinity. What we also prove in this note is the following.

\begin{theorem}\label{thmMAIN1}
Suppose that $u$ is a radially symmetric solution with exactly quadratic growth $\kappa >0$ at infinity found in Theorem \ref{thmMAIN} above. Then we have the following further asymptotic behavior.
\begin{itemize}
  \item[(a)] For $q>3/2$,
\begin{equation*}
\displaystyle\lim_{|x|\to +\infty}\frac{u(x)-\kappa |x|^2}{|x|}=\frac{1}{2}\int_0^{\infty}|x|^2 u^{-q}(x) dx.
\end{equation*}
  \item[(b)] For $q=3/2$,
\begin{equation*}
\lim_{|x|\to +\infty}\frac{u(x)-\kappa |x|^2}{|x| \log(|x|)}=\frac{1}{2\kappa^{3/2}}.
\end{equation*}
  \item[(c)] For $1<q<3/2$,
\begin{equation*}
\lim_{|x|\to +\infty}\frac{u(x)-\kappa |x|^2}{|x|^{4-2q}}=\chi,
\end{equation*}
where
\[
\chi=\frac{1}{2\kappa^q}\left( \frac{1}{3-2q}- \frac{1}{4-2q}+\frac{1}{3(5-2q)}-\frac{1}{3(2-2q)} \right).
\]
\end{itemize}
\end{theorem}

To prove Theorems \ref{thmMAIN} and \ref{thmMAIN1}, we closely follow the argument presented in \cite{G}. Before closing this section, it is worth noting that Theorem \ref{thmMAIN} complements all mentioned results above and hence completes the picture of radially symmetric solutions of \eqref{eqMain}. For clarity, we summary all results above as in Table \ref{tbl}.

\restylefloat{table}
\begin{center}
\begin{table}[H]\label{tbl}
\begin{tabularx}{\linewidth}{ccccc}\hline
 \multicolumn{1}{p{0.16\linewidth}|}{\centering $1<q<3$}         & \multicolumn{1}{p{0.17\linewidth}|}{\centering $q=3$}              & \multicolumn{1}{p{0.17\linewidth}|}{\centering $3 < q < 7$}                & \multicolumn{1}{p{0.18\linewidth}|}{\centering $q=7$}         & \multicolumn{1}{p{0.16\linewidth}}{\centering $q>7$}\\
  \hline \hline
\cellcolor{blue!15}  & \cellcolor{blue!15} &\multicolumn{3}{|c}{\textit{necessary} if $u $ grows linearly} \\
  \hline
 \cellcolor{blue!15} & \cellcolor{blue!15} & \cellcolor{blue!15} & \multicolumn{1}{|c|}{$u $ is \textit{precise}} & \cellcolor{blue!15} \\
  \hline
  \cellcolor{blue!15}  &   \cellcolor{blue!15}  &\multicolumn{3}{|c}{$u$ grows \textit{either linearly or quadratically}} \\
  \hline
\multicolumn{5}{c}{$u(r)$ grows \textit{between linear and quadratic} }  \\
  \hline
  \cellcolor{blue!15}  &   \cellcolor{blue!15}  &  \cellcolor{blue!15}  &\multicolumn{2}{|c}{\textit{$\exists$}! $u $ grows linearly} \\
  \hline
 \cellcolor{blue!15}  &   \cellcolor{blue!15}  &\multicolumn{1}{|c|}{\textit{$\exists$}$u $  linearly } &  \cellcolor{blue!15}  &  \cellcolor{blue!15}  \\
  \hline
  \cellcolor{blue!15}  &\multicolumn{1}{|c|}{\textit{$\exists$} $u \approx r (\log r)^{\frac 14}$ } &  \cellcolor{blue!15}  &  \cellcolor{blue!15}  &  \cellcolor{blue!15}  \\
 \hline
\multicolumn{1}{c|}{\textit{$\exists$}$u \approx r^{4/(1+q)}$}  & \cellcolor{blue!15} & \cellcolor{blue!15} & \cellcolor{blue!15} & \cellcolor{blue!15} \\
  \hline
\multicolumn{5}{c}{\textit{$\exists$} infinitely many $u$ grows \textit{quadratically}}  \\
\hline
\end{tabularx}\medskip 
\caption{Summary of results for radially symmetric solutions $u(r)$ of \eqref{eqMain}.}
\end{table}
\end{center}

\vspace{-20pt} (The first three rows are due to Choi--Xu \cite{ChoiXu}, the next two rows are due to McKenna--Reichel \cite{KR}, the next three rows are due to Guerra \cite{G}, and the last is from Theorem \ref{thmMAIN}.) In one way or another, we know that $q>1$ is necessary and radially symmetric solutions of \eqref{eqMain} must grow linearly up to quadratically. For quadratic growth at infinity, Theorem \ref{thmMAIN} conclude that there are radially symmetric solutions of \eqref{eqMain} which have quadratic growth at infinity for all $q>1$. 

\section{Proof of Theorems \ref{thmMAIN} and \ref{thmMAIN1}}

\subsection{An initial value problem}

The present proof follows the arguments in \cite{G} closely. First, suppose $\beta > 0$, we consider the following initial value problem:
\begin{equation}\label{eqIVP}
\left\{
\begin{split}
\Delta^2 U &= - U^{-q}, \qquad U>0, \qquad r \in (0, R_{\max}(\beta)),\\
U(0)&=1, \quad U'(0)=0, \quad \Delta U(0) = \beta, \quad (\Delta U)'(0)=0,
\end{split}
\right.
\end{equation}
where $[0, R_{\max}(\beta))$ is the maximal interval of existence of solutions. (Such an existence of solutions for \eqref{eqIVP} follows from standard ODE theory. A similar problem in $\mathbb R^2$ and $\mathbb R^3$ for $q=2$ was studied in \cite{GuoWei}.) The following result, indicating the threshold for $\beta$, was obtained in \cite[Proposition 2.1]{G}.

\begin{proposition}\label{propIVP}
Assume that $q>1$ and $\beta>0$. Let $U_\beta$ be the unique local solution of \eqref{eqIVP} above. Then there is a unique $\beta^\star>0$ such that:
\begin{itemize}
  \item [(a)] If $\beta < \beta^\star$ then $R_{\max} (\beta) < \infty$.
  \item [(b)] If $\beta \geqslant \beta^\star$ then $R_{\max} (\beta) = \infty$.
  \item [(c)] If $\beta \geqslant \beta^\star$ then $\lim_{r \to +\infty} \Delta U_\beta (r) \geqslant 0$.
  \item [(d)] We have $\beta = \beta^\star$ if and only if $\lim_{r \to +\infty} \Delta U_\beta (r) = 0$.
\end{itemize}
\end{proposition}

In the rest of our present proof, we set $u=U_\beta$ for some fixed $\beta > \beta^\star$ but arbitrary. Then it suffices to show that $u$ has exactly quadratic growth at infinity. The fact that such a limit at infinity can be arbitrary follows from a suitable scaling of $u$. As a key step toward this end, we shall study asymptotic behavior of $u$ in the next subsection.

\subsection{Asymptotic behavior}

To understand the structure of radially symmetric solutions of \eqref{eqMain}, we transform \eqref{eqMain} into the following system of second order partial differential equations
\begin{equation}\label{eqMainSystem}
\left\{
\begin{split}
\Delta u &= v \quad \text{ in } \mathbb R^3,\\
\Delta v &= - u^{-q} \quad \text{ in } \mathbb R^3.
\end{split}
\right.
\end{equation}
To study the asymptotic behavior of \eqref{eqMainSystem}, we follow the ideas in \cite{HV}. First, by the Emden--Fowler transformation we set
\begin{equation}\label{eqVariableChange}
x(t) = \frac{ru'}u, \quad y(t) = \frac{rv'}v, \quad z(t)=\frac{r^2 v}u, \quad w(t)=\frac{r^2 u^{-q}}v, \quad t =\log r.
\end{equation}
Then the system \eqref{eqMainSystem} is transformed into a $4$-dimensional quadratic system of the form
\begin{equation}\label{eqSystem}
\left\{
\begin{split}
x' &=x(-1-x)+z,\\
y &'=y(-1-y)-w,\\
z &'=z(2-x+y),\\
w'&=w(2-qx-y),
\end{split}
\right.
\end{equation}
where $'=d/dt$. As indicated in \cite{G}, the critical points of \eqref{eqSystem} are
\[
\begin{split}
p_0 & = (0,0,0,0), \quad p_1 =(1,-1,2,0), \quad p_2=(2,0,6,0),\\
p_3 &= (a,a-2,a(a+1),(2-a)(a-1)), \quad p_4 = (0,2,0,-6), \quad p_5  = (0,-1,0,0),\\
p_6 &=(-1,0,0,0),\quad  p_7 = (-1,-1,0,0), \quad p_8 = (-1,q+2, 0,-(q+2)(q+3)),
\end{split}
\]
where $a=4/(q+1)$.

Thanks to Proposition \ref{propIVP}, the solution $u$ of \eqref{eqIVP} exists for all time $t$; hence we can denote
\[
\gamma = \lim_{r \to +\infty} v(r),
\]
where, as always, we set $v=\Delta u$. By Proposition \ref{propIVP}(c, d) we know that $\gamma >0$. Since $-\Delta v = u^{-q}$ with $v'(0)=0$, we have
\begin{equation}\label{eqV(r)=}
v(r) = v(0) - \int_0^r t u^{-q} dt + \frac 1r \int_0^r t^2 u^{-q} dt.
\end{equation}
Since $\gamma >0$, there exist two positive constant $a$ and $b$ such that $u(r) \geqslant ar^2 + b$ holds for all $r>0$. From this, for each $q>1$, we find that $t^2 u^{-q} \to 0$ as $t \to +\infty$. Thanks to the l'H\^opital rule, we can pass \eqref{eqV(r)=} to the limit as $r \to +\infty$ to obtain
\begin{equation}\label{gamma=}
v(0)= \gamma +  \int_0^{+\infty} t u^{-q} dt < \infty
\end{equation}
while on the other hand we get
\begin{equation}\label{eq1V(r)=}
v(r) = \gamma+ \int_r^{+\infty} t u^{-q} dt + \frac 1r \int_0^r t^2 u^{-q} dt.
\end{equation}
Now, an easy computation leads us to
\[
\frac{v(r)}{r v'(r)} = - \frac{r \int_r^{+\infty} t u^{-q} dt + \gamma r}{\int_0^r t^2 u^{-q} dt}-1.
\]

\noindent\textbf{Claim 1}. There holds $v(r)/(r v'(r)) \to -\infty$ as $r \to +\infty$. In other words, $y(r) \to 0$ as $r \to +\infty$.

\begin{proof}[Proof of Claim 1]
Depending on the value of $q$, there are two possible cases.

\noindent\underline{Case 1}. Suppose $q>3/2$, then we immediately see that the integral $\int_0^{+\infty} t^2 u^{-q} dt$ converges. Hence the claim holds since $\gamma >0$.

\noindent\underline{Case 2}. In this scenario, there holds $1<q \leqslant 3/2$. Using the l'H\^opital rule, we arrive at
\[
\lim_{r \to +\infty} \frac{v(r)}{r v'(r)} = - \lim_{r \to +\infty}  \frac{\int_r^{+\infty} t u^{-q} dt  + \gamma r}{ r^2 u^{-q} } = -\infty,
\]
thanks to $u^{-q} / r \to +\infty$ as $r \to +\infty$.
\end{proof}

Now using the relation $\Delta u = v$ with $u(0)=1$ and $u'(0)=0$, in a similar fashion of \eqref{eqV(r)=} we obtain
\begin{equation}\label{eqU(r)=}
u(r) = 1 + \int_0^r t v dt - \frac 1r \int_0^r t^2 v dt.
\end{equation}
From this, we obtain
\begin{equation}\label{eqU(r)/R^2V(t)=}
\frac{u(r)}{r^2 v(r)} = \frac{1}{r^2 v(r)} + \frac{1}{r^2 v(r)} \int_0^r t v dt - \frac 1{r^3 v(r)} \int_0^r t^2 v dt.
\end{equation}

\noindent\textbf{Claim 2}. There holds $u(r)/(r^2 v(r)) \to 1/6$ as $r \to +\infty$. In other words, $z(r) \to 6$ as $r \to +\infty$.

\begin{proof}[Proof of Claim 2]
We first use the l'H\^opital rule to estimate the last two terms on the right hand side of \eqref{eqU(r)/R^2V(t)=}. For the middle term, we clearly have
\[
\lim_{r \to +\infty} \frac 1{r^2 v(r)} \int_0^r t v(t) dt =  \lim_{r \to +\infty}  \frac{1}{ 2 + r v'(r) / v(r)} = \frac 12,
\]
thanks to Claim 1. For the last term, we get
\[
\lim_{r \to +\infty} \frac 1{r^3 v(r)} \int_0^2 t v(t) dt =  \lim_{r \to +\infty}  \frac{1}{ 3 + r v'(r) / v(r)} = \frac 13.
\]
We are now in a position to estimate $u(r)/(r^2 v(r))$ when $r$ is large. Clearly, by \eqref{eqU(r)/R^2V(t)=} we know that $\lim_{r \to +\infty} u(r)/(r^2 v(r)) = 1/6$ since $r^2 v(r) \to +\infty$ as $r \to +\infty$.
\end{proof}

\noindent\textbf{Claim 3}. There holds $u(r)/(r u'(r)) \to 1/2$ as $r \to +\infty$. In other words, $x(r) \to 2$ as $r \to +\infty$.

\begin{proof}[Proof of Claim 3]
Using \eqref{eqU(r)=}, we obtain $u'(r) = r^{-2} \int_0^r t^2 v dt$ which yields
\begin{equation}\label{eqU(r)/RU'(t)=}
\frac{u(r)}{r u'(r)} = \frac{ r}{\int_0^r t^2 v(t) dt} + \frac{r \int_0^r t v(t) dt }{ \int_0^r t^2 v(t) dt} -1.
\end{equation}
The l'H\^opital rule applied to the first term on the right hand side of \eqref{eqU(r)/RU'(t)=} gives
\[
\lim_{r \to +\infty}  \frac{ r}{\int_0^r t^2 v(t) dt} = 0
\]
while for the second term we know that
\[
\lim_{r \to +\infty} \frac{r \int_0^r t v(t) dt }{ \int_0^r t^2 v(t) dt} = \frac 32.
\]
From this we obtain the desired limit.
\end{proof}

From Claims 1, 2, and 3 we see that the solutions $(x,y,z,w)$ corresponding to the radially symmetric solutions with 
quadratic growth are attracted to the point $p_2:=(2,0,6,0)$ at infinity. Therefore, the asymptotic behavior is obtained by analyzing
the asymptotic behavior of solutions about $(2,0,6,0)$. 

For $q>1$, we first obtain the linearization of \eqref{eqSystem} at $(x,y,z,w)$ given by the following matrix
\[
\begin{pmatrix}
-2x-1  & 0 & 1 & 0\\
0 & -2y-1  & 0 & 1\\
-z  & z & -x+2 & 0\\
-qw  & -w & 0 & -qx-y-2
\end{pmatrix}.
\]
At $p_2$, this matrix becomes
\[
\begin{pmatrix}
-5  & 0 & 1 & 0\\
0 & -1  & 0 & 1\\
-6  & 6 & 0 & 0\\
0  & 0 & 0 & -2q -2
\end{pmatrix}
\] 
which has the following eigenvalues: $\lambda_1 =-1$, $\lambda_2 =-2$, $\lambda_3 =-3$, and $\lambda_4 =2-2q$. Since these eigenvalues are non-zero whenever $q>1$, we conclude that there exists a constant $c_q \ne 0$ such that the following asymptotic behavior occurs: For $q>3/2$
\begin{equation}\label{Case q>3/2}
\frac{ru'(r)}{u(r)}=2+ c_q e^{-t} +o(e^{-t}) 
\end{equation}
as $t\to +\infty$ while for $q= 3/2$
\begin{equation}\label{Case q=3/2}
\frac{ru'(r)}{u(r)}=2+ c_q t e^{-t} +o(t e^{-t}) 
\end{equation}
as $t\to +\infty$ due to $\lambda_1 = \lambda_4$, and for $1<q<3/2$
\begin{equation}\label{Case q<3/2}
\frac{ru'(r)}{u(r)}=2+ c_q e^{-(2q-2)t} +o(e^{-(2q-2)t}) 
\end{equation}
as $t\to +\infty$. 

\subsection{Quadratic growth at infinity can be arbitrary}

We establish in this subsection the fact that if there is some radially symmetric solution $u$ of \eqref{eqMain} having 
$$ \lim_{|x|\to +\infty} \frac{u(r)}{r^2}=\kappa$$ 
for some $\kappa > 0$, then given any $\varpi>0$, there exists a radially symmetric solution $v$ of \eqref{eqMain} such that
$$ \lim_{|x|\to +\infty} \frac{v(r)}{r^2}=\varpi.$$ 
To see this, we first set 
\[
v(r) =\Big(\frac \varpi\kappa \Big)^\delta u\Big(\big(\frac \varpi\kappa \big)^\alpha r\Big).
\]
Then it is elementary to see that $v$ solves \eqref{eqMain} if $(1+q)\delta + 4\alpha=0$. To fulfill the limit $\lim_{|x|\to +\infty} v(r)/r^2=\varpi$, it also requires $\delta+2\alpha =1$. Resolving these conditions for $\delta$ and $\alpha$, we conclude that $\delta = -2/(q-1)$ and $\alpha = (q+1)/(2(q-1))$ which are obviously well-defined for all $q>1$.

\subsection{Proof of Theorem \ref{thmMAIN}}

Combining \eqref{eq1V(r)=} and \eqref{eqU(r)=}, we have the following representation
\begin{equation}\label{eqdenote u(r)=}
u(r)=\frac{r}{2}\int_0^{r}t^2 u^{-q}dt -\frac{1}{2}\int_0^r t^3 u^{-q}dt +\frac{1}{6r} \int_0^r t^4 u^{-q}dt
+\frac{r^2}{6} \int_r^{+\infty} t u^{-q}dt +\frac{\gamma \, r^2}{6} +1
\end{equation}
which is similar to \cite[Eq. (5.1)]{ChoiXu}. We note that the representation \eqref{eqdenote u(r)=} is valid for all $q>1$. In general, the term $\int_r^{+\infty} t u^{-q}dt$ may not be well-defined for $q<2$ if $u$ solves \eqref{eqMain}.

Note that with $r=e^t$ we obtain $h'(t)= ru'(r)/u(r)$ where we set $h(t):=\log u(r)$. Therefore, by using \eqref{Case q>3/2}--\eqref{Case q<3/2}, we obtain
\[
h'(t)= 
\begin{cases}
2+ c_q e^{-t} +o(e^{-t}) & \text{ if } q>3/2\\
2+ c_q t e^{-t} +o(t e^{-t})  & \text{ if } q=3/2\\
2+ c_q e^{-(2q-2)t} +o(e^{-(2q-2)t})  & \text{ if } 1<q<3/2
\end{cases}
\]
as $t \to +\infty$. Integrating both sides gives
\[
\frac {u(r)}{r^2} = 
\begin{cases}
u(1) \exp \int_0^t \big(   c_q e^{-s} +o(e^{-s}) \big) ds & \text{ if } q>3/2\\
u(1) \exp \int_0^t \big(   c_q s e^{-s} +o(t e^{-s})  \big) ds & \text{ if } q=3/2\\
u(1) \exp \int_0^t \big(   c_q e^{-(2q-2)s} +o(e^{-(2q-2)s})  \big) ds & \text{ if } 1<q<3/2
\end{cases}
\]
From this, it is easy to see that the following limit $ \lim_{|x|\to +\infty} u(r)/r^2=\kappa$ exists for some $\kappa > 0$. Now using \eqref{eqdenote u(r)=} and \eqref{gamma=}, we further obtain 
\begin{equation}\label{kappa=}
\kappa =\frac{\gamma}{6}=\frac{1}{6}\left( (\Delta u)(0) -\int_0^{+\infty}t u^{-q}(t) dt \right)
\end{equation} 
as claimed. The fact that $\kappa$ can be arbitrary follows from the preceding subsection by scaling $u$. 

Finally, since one can freely choose $\beta > \beta^\star$, we conclude the existence of infinitely many radially symmetric solutions \eqref{eqMain} of exactly given quadratic growth at infinity for any $q>1$. To realize this fact, one first pick $\beta_1 \ne \beta_2 > \beta^\star$ and follow the procedure above to select $u_{\beta_i}$ with growth $\kappa_i$ at infinity, that is $ \lim_{|x|\to +\infty} u_{\beta_i} (r)/r^2=\kappa_i$ for $i=1,2$. Then we define
\[
w_{\beta_2} =\begin{cases}
u_{\beta_2}& \text{ if } \quad \kappa_1 = \kappa_2, \\
 \Big(\dfrac {\kappa_1} {\kappa_2} \Big)^{-\frac 2{q-1}} u_{\beta_2} \Big(\big(\dfrac {\kappa_1} {\kappa_2} \big)^\frac{q+1}{2(q-1)} r\Big) &\text{ if } \quad \kappa_1 \ne \kappa_2.
\end{cases}
\]
Note that $u_{\beta_i}$ are distinct because $\beta_1 \ne \beta_2$. Moreover, in the case $\kappa_1 \ne \kappa_2$ there holds $w_{\beta_2}(0)= (\kappa_1 /\kappa_2)^{-2/(q-1)} \ne 1$; hence $u_{\beta_1} \not\equiv w_{\beta_2}$. In addition, it is easy to see that $\lim_{|x|\to +\infty} w_{\beta_2}(r)/r^2 = \kappa_1$. Therefore, given $\varpi>0$ if we scale $u_{\beta_1}$ and $w_{\beta_2}$ to get
\[
v_1(r) =\Big(\frac \varpi{\kappa_1} \Big)^{-\frac 2{q-1}} u_{\beta_1} \Big(\big(\frac \varpi{\kappa_1} \big)^\frac{q+1}{2(q-1)} r\Big)
\]
and
\[
v_2(r) =\Big(\frac \varpi{\kappa_1} \Big)^{-\frac 2{q-1}} w_{\beta_2} \Big(\big(\frac \varpi{\kappa_1} \big)^\frac{q+1}{2(q-1)} r\Big),
\]
it is immediate to see that $v_1 \not\equiv v_2$ and that 
\[
\lim_{|x|\to +\infty} \frac{v_1(r)}{r^2}= \lim_{|x|\to +\infty} \frac{v_2(r)}{r^2} = \varpi.
\]
The proof is complete.

\subsection{Proof of Theorem \ref{thmMAIN1}}

Let $u$ be a  solution of \eqref{eqMain} constructed as above which has quadratic growth $\kappa$ at infinity. Using the quadratic growth formula \eqref{kappa=} and the presentation \eqref{eqdenote u(r)=}, we know that the constructed solution $u$ also fulfills following presentation
\begin{equation}\label{eq u(r)-kappa r^2=}
u(r) - \kappa r^2 =\frac{r}{2}\int_0^{r}t^2 u^{-q}dt -\frac{1}{2}\int_0^r t^3 u^{-q}dt +\frac{1}{6r} \int_0^r t^4 u^{-q}dt
+\frac{r^2}{6} \int_r^{+\infty} t u^{-q}dt  +u(0),
\end{equation}
and this presentation is also valid for all $q>1$, compared with \cite[Eq. (5.1)]{ChoiXu}. 

Then we make use of \eqref{Case q>3/2}--\eqref{Case q<3/2} plus the l'H\^opital rule to conclude the theorem. For example, when $q>3/2$, there holds $(u(r) - \kappa r^2)/r \to \int_0^{+\infty}t^2 u^{-q}dt$ due to the contribution of the first integral in \eqref{eq u(r)-kappa r^2=} since $r^3 u^{-q}(r) \to 0$ as $r \to +\infty$. When $q=3/2$, $(u(r) - \kappa r^2)/(r \log r) \to 1/(2\kappa^{3/2})$ due to the contribution of the second integral in \eqref{eq u(r)-kappa r^2=} while in the case $q<3/2$, $(u(r) - \kappa r^2)/r^{4-2q} \to \chi$ due to the contribution of all four integrals in \eqref{eq u(r)-kappa r^2=}. This completes our proof of Theorem \ref{thmMAIN1}.

\section*{Acknowledgments}

We thank Prof. Ignacio Guerra for useful discussion about his paper \cite{G}. The second author was partially supported by Vietnam National University at Hanoi through VNU Scientist Links.


\end{document}